# Teaching mathematics with a different philosophy

# Part 1: Formal mathematics as biased metaphysics

C. K. Raju
School of Mathematical Sciences
Universiti Sains Malaysia
11700, Penang, Malaysia
ckr@ckraju.net

**Abstract:** We report on a pedagogical experiment to make mathematics easy by changing its philosophy. The Western philosophy of math originated in religious beliefs about mathesis, cursed by the church. Later, mathematics was "reinterpreted", in a theologically-correct way, using the myth of "Euclid" and his deductive proofs. The *fact* of the empirical proofs in the *Elements,* was, however, contrary to this myth. The discrepancy was resolved by Hilbert and Russell who rejected empirical proofs as unsound, reducing all mathematics to metaphysics. We explain why that formalist metaphysics is anti-utilitarian and culturally biased, *not* universal. Historically, most school-level math originated in the non-West with a practical epistemology, but was absorbed in the West after superimposing on it an incompatible Western metaphysics, still used to teach it. This has made mathematics needlessly complex. Accordingly, math can be made easy and more universal by reverting to a more practical epistemology.

## *1. Introduction*

Learning difficulties with math are widespread. This paper reports on a study which has (a) identified the origin of the learning difficulties with math in a novel way, (b) proposed an alternative, and (c) tested the feasibility of teaching that alternative for the case of the calculus course currently taught at the 12$^{th}$ std and first year university level.

## *2. The analysis*

To remedy an ailment, the first step is to understand its causes. Why, then, is math difficult? Learning difficulties with mathematics are offhand put down to sundry causes—lack of good teachers, lack of "mathematical aptitude", etc. No doubt an inspired teacher can make any subject interesting; however, such "causes" do not explain why so many school students fear math, and not geography, for example.

The new answer locates the cause in the nature of the subject. The Western philosophy of math was religiously oriented; hence theological complexities got intertwined with present-day math. This has made math difficult.

The idea of math as religiously oriented may seem surprising. However, the very word "mathematics" derives from "mathesis", which means recollection of knowledge *from previous lives*. In Plato's *Meno*,[1] Socrates demonstrates mathesis by questioning an uneducated slave boy to elicit his innate knowledge of mathematics. This, concludes Socrates, is proof that the boy has an immortal soul which got that innate knowledge from previous lives. Why math, and not some other form of knowledge? Because mathematics was believed to incorporate *eternal* truths which most readily aroused the eternal (immortal) soul, as Plato stated.[2] Arousing the soul was thought to make people virtuous. Hence, in his *Republic*, Plato advocated the teaching of mathematics to make people virtuous; he explicitly rejected its practical applications as unimportant. This idea of math as good for the soul persisted for a thousand years: Proclus (5th c.) too explained that math leads to "the blessed life".[3]

This connection of mathematics to beliefs about the soul drew it into a religious conflict. The Christian church initially had similar beliefs about the soul and its previous lives.[4] However, after marrying the state, in the 4th c., the church changed its notion of soul, to suit the needs of state power. Philosophers who persisted with the earlier beliefs about the soul (underlying mathesis) were violently attacked. Their key resources (libraries and temples[5]) were all destroyed.[6] Those mathematicians who continued to resist, like Hypatia, were lynched. Her successor, Proclus was declared a heretic for his belief that the cosmos must be eternal hence uncreated[7] (since the truths of math are eternal). In 552 CE, the church cursed the belief in past lives underlying mathesis; those curses are known today as the "anathemas against pre-existence".[8] Philosophy and mathematics were banished from Christendom.

Mathematics and philosophy were eventually accepted back during the Crusades, but this entailed an attempt to resolve the theological conflict by "reinterpreting" the philosophy of math in a theologically-correct way. This reinterpretation eliminated mathesis from math, making it soul-less, so to say. It was claimed that math was *only* about reasoning and proof, thus aligning the philosophy of math with the post-Crusade Christian theology of reason.[9]

This reinterpretation of mathematics is closely tied to a false history[10, 11] which attributes geometry to "Euclid", and claims that his concern was not mathesis but "irrefragable demonstration". What do we really *know* about this "Euclid"? "Nothing" as a Western authority on Greek mathematics, the late David Fowler, publicly admitted.[12] As for the book, *Elements*, the available papyri fragments[13] show that it was *not* standardised even seven centuries after the date of "Euclid" who had supposedly standardised it. The name "Euclid" is *not* mentioned in *any* Greek manuscript of the *Elements,* which *all* attribute it to Theon (4th c.), or say it is based on his lectures.[14] Even the Greek commentaries do not mention "Euclid", but speak *anonymously* of "the author of the Elements"[15] (perhaps because the author was a woman, Hypatia, Theon's daughter). Hypatia's viewpoint, like that of her successor Proclus, was still that of mathesis. So, attributing the *Elements* to "Euclid" is critical to the post-Crusade philosophy of mathematics which delinked math from mathesis.

Crusading historians had another strong motive to falsely attribute texts to early[16] Greeks.[17] At this time of intense religious fanaticism, the church found it difficult to abandon its earlier tradition[18] of burning heretical books, and instead translate books captured from the enemy. This policy change was justified by the story that all Arabic scientific books were written by early Greeks, and hence were the rightful inheritance of Europe which Arabs had only preserved. (Later day racist and colonial historians built upon this distorted history which persists to this day in school texts.[19]) The name "Euclid" is found only in Latin translations from the Arabic, where it may have arisen as a translation mistake: a misreading of the Arabic "ucli des" (meaning "key to geometry") as Uclides, the name of a Greek.[20]

The present-day idea of mathematics as proof originates from "Euclid's" purported philosophy of "irrefragable demonstration". All the evidence for that philosophy is an isolated passage in a late manuscript[21] of Proclus' *Commentary.* That passage is spurious. Thus, it claims that Archimedes cited "Euclid". Such a citation (of the *Elements*, not "Euclid") is indeed found in a 15th c. manuscript of the *Sphere and Cylinder*,[22] somehow attributed to Archimedes. However, citing texts in this modern manner of Christian scribes was not the custom in Archimedes' time. Further, there is only one citation, though are many other places in that work where the *Elements* could have been cited, but is not. Hence, that isolated "Archimedes citation" is recognized[23] as spurious. Since the author of the "Proclus passage" knew of the spurious "Archimedes citation", the "Proclus passage" must itself be a spurious interpolation from even later. As such, there is no serious basis for the belief that Greeks ever regarded mathematics as concerned with deductive proof rather than mathesis.

Since so much has been anachronistically read into those two words "irrefragable demonstration", it should also be mentioned that the (spurious) passage itself claims that the "irrefragable demonstration" was based on "causes and [astrological?] signs"! In any case, the myth persisted that "Euclid" intended a special metaphysical notion of (deductive) proof, and that myth is used to justify the present-day belief that such proofs (rather than calculation) are the key concern of mathematics.

That myth is contrary to facts, for the *Elements did* make essential use of empirical (or physical) means of proof in the 1$^{st}$ and the 4$^{th}$ (side-angle-side) proposition, to mention two well-known[24] examples. If empirical proof is allowed in one place, why not in another? Why not prove the 47$^{th}$ proposition ("Pythagorean theorem") directly, by empirical means? That greatly simplifies the proof: it requires only one step instead of the 46 intermediate propositions actually used in the *Elements*. (That was how it was proved, empirically, in one step, in Indian texts, like the *Yuktibhasa*, where it is the first proposition, not the last.) So, this myth (that the *Elements* concerns metaphysical proof) makes the *Elements* either trivial or inconsistent. Proclus could explain that "proofs must vary...with the kinds of being",[25] for he regarded the *Elements* as a *religious* book, concerning mathesis. But the inconsistent means of proof cannot be explained on the post-Crusade reinterpretation of the *Elements* (and mathematics) as *primarily* concerned with metaphysical proof.

The present-day philosophy of formal math arose directly from the attempts by Hilbert and Russell to reconcile (a) the *myth* of "Euclid" and his deductive proofs with (b) the *fact* of the empirical proofs in the *Elements*. The myth prevailed over the facts. In their respective tracts[26] on the foundations of geometry both tried to "save the story" of "Euclid". They rejected those empirical proofs in the *Elements* as unsound, since incompatible with the (purported) intentions of (the mythical) "Euclid"! Since the side-angle-side (SAS) *theorem* is essential, they changed it into a *postulate.* (And that is how geometry has been taught in schools[27] since the 1970's, using the SAS postulate.) Formalism developed from that analysis: on the same lines, Hilbert and Russell made *all* math completely formal and metaphysical as it is today.

Formalism made even elementary math difficult, although Russell regarded that[28] as its chief advantage! (Recall that the materialist Epicureans had criticised geometry saying that its theorems were obvious even to an ass; so, non-obviousness justifies the philosophy that those theorems needed proofs.) The postulates used by Hilbert to formalise "Euclidean" geometry were similarly *intended* to make geometry difficult. Thus, "Euclidean" geometry may also be formalised using Birkhoff's metric postulates.[29] But, that again trivialises the *Elements* (*if* it is regarded as primarily concerned with proof), for proofs are again simplified. But if the *Elements* is trivial*,* then formalism based on it would itself seem ridiculous. Hence, Hilbert proposed synthetic geometry which rejects the notion of length. This makes elementary geometry much harder but justifies, in a general way, why the theorems in the *Elements* needed complex proofs. However, it does not fit the particulars: synthetic geometry is *not*[30] a valid interpretation of the *Elements*. Thus, Hilbert's synthetic term "congruence" is *not* a valid substitute for the original word "equality" used in the *Elements*. Proposition 35 of the *Elements* is about the equality of the *areas* of *incongruent* parallelograms (on the same base and between the same parallels). It is still possible to "save the story" by claiming that "Euclid" defined area, but not length (and hence needed lengthy proofs). But, this is manifestly far-fetched. Moreover, it is well known in the philosophy of science that *any* theory/story can be saved from *any* facts by sufficiently many hypotheses: one lie can always be defended by telling a thousand more. Hence, such a long chain of hypotheses is not worth discussing further; it is simpler to reject the story.

Setting aside the question of its historical origins in the *Elements*, one can enquire about the *utility* of formalism. Does formal math *suit* applications to science and technology? The idea of math as metaphysics (which is superior to physics) naturally makes it *less* suited for applications to science and technology which concern the physical and empirical. For example, a discontinuous function is *not* differentiable on elementary mathematical analysis, but *is* differentiable on the Schwartz theory of distributions. That is, formal mathematics being metaphysical, a discontinuous function is differentiable or not, as one likes! So what should one believe about the partial differential equations of physics? Do they or do they not apply to[31] observed discontinuities (like shock waves) or the related "singularities"? Actually, both definitions of the derivative fail under these circumstances (the one because discontinuous functions are not differentiable, and the other because the equations of physics are nonlinear and Schwartz distributions cannot be multiplied). So should we reject them both? Clearly, this believe-what-you-like approach is unsuited to scientific theory which is expected to be definite and refutable. The only other option is to allow empirical inputs into mathematics. Selecting a definition of the derivative based on empirical inputs would explicitly make mathematics an adjunct physical theory, as in this author's philosophy. Empirical inputs always suit science and technology. The alternative of banking on mathematical authority is not reliable, as the case of Riemann shows. (He mistakenly[32] assumed that entropy is constant across a shock.)

Does formalism, then, provide a universal metaphysics? Now, it is an elementary matter of commonsense that metaphysics can never be universal. However, the case of 2+2=4 is often naively cited as "proof" of the universality of mathematics. This is naïve because the *practical* notion of 2 which derives as an abstraction from the *empirical* observation of 2 dogs, 2 stones etc. has nothing whatsoever to do with *formal* mathematics. It is perfectly possible to have a formal theory[33] in which 2+2=5, say, just as 2 male rabbits and 2 female rabbits may make any number of rabbits over a period of time. Likewise, the circuits on a computer chip routinely implement an arithmetic in which 1+1=0 (exclusive disjunction), or 1+1=1 (inclusive disjunction). Thus, formally, it is necessary to specify that the symbols 2 , +, and 4 relate to Peano's postulates. Trying to specify this brings in the metaphysics of infinity—a real computer (with finite memory, not a Turing machine) can *never* implement Peano arithmetic,[34] because the notion of a natural number cannot be finitely specified. Thus, formalism does *not* provide a universal metaphysics. However, the philosophy of mathematics as metaphysics, combined with the myth of mathematics as universal truth, helped to promote a particular brand of metaphysics as universal.

This is problematic because while formal mathematics is no longer *explicitly* religious, like mathesis, its metaphysics remains religiously biased. On post-Crusade Christian rational theology, it was thought that God is bound by logic (cannot create an illogical world) but is free to create empirical facts of his choice. Hence, Western theologians came to believe that logic (which binds God) is "stronger" than empirical facts (which do not bind God). This theological belief is exactly mirrored in the present-day Wittgenstein-Tarski semantics of possible worlds in which logical truths are regarded as *necessary* truths (true in *all* possible worlds) unlike empirical truths, which are regarded as weaker *contingent* truths (true in *some* possible worlds). The only difference is that instead of speaking of "possible worlds which God could create", we speak of "possible worlds according to Wittgenstein"! That is, the metaphysics of formal math is aligned to post-Crusade Western theology which regarded metaphysics as *more* reliable than physics.

In sharp contrast, *all* Indian systems of philosophy, without any exception, accept the empirical (*pratyaksa*) as the first means of proof (*pramana*) while the Lokayata *reject* inference/deduction as unreliable. So, Indian philosophy considered empirical proof as more reliable than logical inference.

Thus, the contrary idea of metaphysical proof as "stronger" than empirical proof would lead at one stroke to the rejection of *all* Indian systems of philosophy. This illustrates how the metaphysics of formal math is not universal but is biased against other systems of philosophy.

Now, deductive inference is based on logic, but *which* logic? Deductive proof lacks certainty unless we can answer this question with certainty. Russell thought, like Kant,[35] that logic is unique and comes from Aristotle. However, one could take instead Buddhist or Jain logic,[36] or quantum logic,[37] or the logic of natural language, none of which is 2-valued. The theorems that can be inferred from a given set of postulates will naturally vary with the logic used: for example, all proofs by contradiction would fail with Buddhist logic. One would no longer be able to prove the existence of a Lebesgue non-measurable set, for example. *This conclusively establishes that the metaphysics of formal math is religiously biased*, for the theorems of formal mathematics vary with religious beliefs. Furthermore, the metaphysics of formal math has no other basis apart from Western culture: it can hardly be supported on the empirical grounds it rejects as inferior!

The religious bias also applies to the *postulates*. In principle, a formal theory could begin with any postulates. However, in practice, those postulates are decided by authoritative mathematicians in the West, as in Hilbert's synthetic geometry. The calculus, as taught to millions of school students today, is based on the notion of limits and the continuum. As noted by Naquib al-Attas, the idea of an infinitely divisible continuum is contrary to the beliefs of Islamic thinkers like al-Ghazali and al-Ashari who believed in atomism. (In fact, the calculus originated in India with similar atomistic beliefs: that the subdivisions of a circle must stop when they reach atomic proportions.) This does not affect *any* practical application: all practical application of calculus today can be done using computers which use floating point numbers which are "atomistic", being finite. Neglecting small numbers is *not* necessarily erroneous, since it is not very different from neglecting infinitesimals in a non-Archimedean field, and can be similarly formalised, since the (formal) notion of infinitesimal is not God-given but is a matter of definition; but students are never told this. They are told that any real calculation done on a computer is forever erroneous, and the only right way to do arithmetic is by using the metaphysics of infinity built into Peano's postulates or the postulates of set theory. Likewise, they are taught that the only right way to do calculus is to use limits. Thus, school students get indoctrinated with the Western theological biases about infinity built into the notions of formal real numbers and limits, which notions are of nil practical value for science and engineering which require real calculation.

In contrast to this close linkage of mathematics to theology in the West, most school math (arithmetic, algebra, trigonometry, calculus) actually originated in the non-West for practical purposes. From "Arabic numerals" (arithmetic algorithms) to trigonometry and calculus, this math was imported by the West for the practical advantages it offered (to commerce, astronomy, and navigation).[38] This practically-oriented non-Western mathematics (actually *ganita*, or *hisab*), had nothing to do with religious beliefs such as mathesis. However, because of its different epistemology[39] it posed difficulties for the theologically-laced notion of mathematics in the West. For example, from the *sulba sutra* to Aryabhata to the *Yuktibhasa*, Indian mathematics freely used empirical means of proof. Obviously, an empirical proof will not in any way diminish the *practical* value of mathematics. However, trying to force-fit this practical, non-Western math into Western religiously-biased ideas about math as metaphysical made the simplest math enormously complicated. Colonisation globalised Western ideas of mathematics, and they continued to be taught today without critical re-examination. On the principle that phylogeny is ontogeny, the actual teaching of mathematics replays in the classroom the historical Western difficulties in absorbing non-Western math.[40] That is what makes math difficult today.

The preceding considerations are illustrated by an example—the calculus—in part 2 of this article.